# Color fixation, color identity and the Four Color Theorem

Asbjørn Brændeland

I argue that *there is no 4-chromatic planar graph with a joinable pair of color identical vertices*, i.e., given a 4-chromatic planar graph $G$ and a pair of vertices $\{u, v\}$ in $G$, if the color of $u$ equals the color of $v$ in every 4-coloring of $G$, then there is no planar supergraph of $G$ where $u$ and $v$ are adjacent.[1]

This is equivalent to the Four Color Theorem, which says that *no more than four colors are required to color a planar map so that no two adjacent regions on the map have the same color*, or, equivalently, that *every planar graph is 4-colorable*.

**Definition 1**. Two vertices $u$ and $v$ in a $k$-chromatic graph $G$ are **color identical** if and only if the color of $u$ equals the color of $v$ in every $k$-coloring of $G$.

**Definition 2**. Color identity is a result of *coloring constraints*. Let $S$ and $T$ be subgraphs of a $k$-chromatic graph such that every vertex in $S$ is adjacent to every vertex in $T$. Then $S$ *constrains* $T$ to the colors $S$ does not have, and vice versa, and if $S$ and $T$ together requires $k$ colors, then $S$ **color fixes** $T$, and vice versa. In particular, if $S$ is an odd cycle and $T$ is a single vertex $v$ in a 4-chromatic graph, and another vertex $u$ is also adjacent to every vertex in $S$, then $S$ fixes $u$ and $v$ to the same color, and we have a color identical pair.

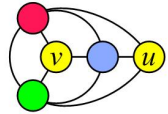

Figure 1.

$S$ may be **color fixed as a whole**, i.e. $S$ may be fixed to $j$ colors in virtue of (**a**) being $j$-critical, and (**b**) having all its vertices constrained to the same set $C$ of $j$ colors, or, to be precise, having at least one vertex constrained to all of $C$ and all vertices constrained to some of $C$. Otherwise, the vertices of $S$ are **color fixed independently**.

If $S$ color fixes $T$ then $\{V(S), V(T)\}$ is a **color fixation pair**, and if $S$ is $j$-critical and $T$ is $(k - j)$-critical, $j < k$, then $\{S, T\}$ is a **color fixation embrace**.

The distinction between *pairs* and *embraces* reflects the fact that not all color fixation pairs induce color fixation embraces, as demonstrated by the 3-chromatic graph in Figure 2. However, as argued below[2], every color fixation pair in a 4-chromatic graph induces a color fixation embrace.

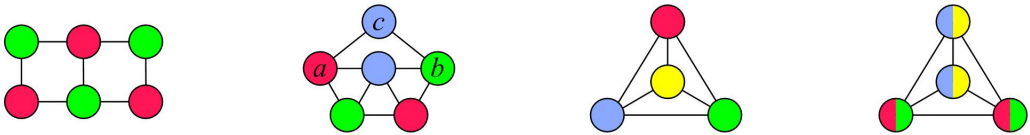

Figure 2. In the 2-chromatic graph to the left, every edge is a color fixation embrace of two vertices. In the 3-chromatic graph, {{a, b}, {c}} is a color fixation pair, but since *a* and *b* are not adjacent, the pair does not induce a color fixation embrace. The two 4-chromatic graphs to the right are both color fixation embraces, of a vertex and a triangle (due to symmetry, we have in fact four such pairs), and of two edges, respectively.

Color fixations are by their nature *local*, but they can *propagate* through chains of color fixation pairs and become *global* (or *less local*), in the sense that two *distant* color fixations can be in the same color fixation chain. What is seen as the primary color fixator in such a chain, is a matter of choice.

**Definition 3**. In the description of the relation between two color fixed subgraphs *A* and *B* of a graph *G*, we can use a ***coloring reference graph***, or ***CR graph***, i.e. a subgraph *R* of *G* that (**i**) we can color fix by choice, without committing to a particular coloring, and (**ii**) whose color fixation explains the color fixation of both *A* and *B*—which is to say that *R* color fixes both *A* and *B* directly or indirectly.

Only a (*k*–1)-complete graph satisfies condition (**i**), since it is the only graph that always receives exactly *k* – 1 colors, and the only graph that receives one distinct color per vertex. Thus every *edge* in a 3-chromatic graph and every *triangle* in a 4-chromatic graph satisfies (**i**).[3]

If *G* has no (*k*–1)-complete subgraph that satisfies condition (**ii**) then *A* and *B* are not both global.

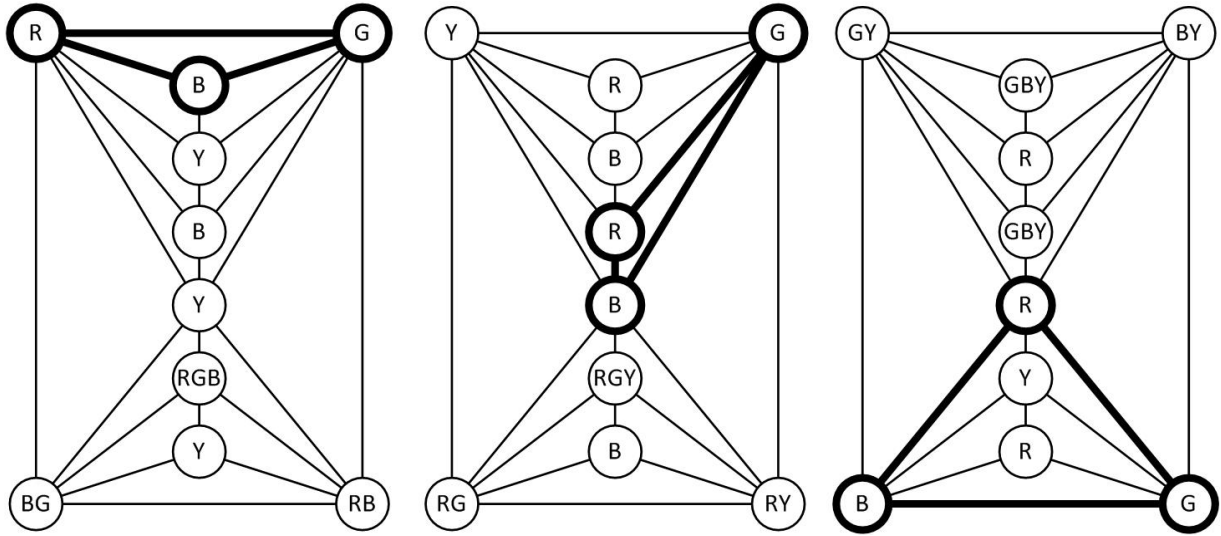

Figure 2. We have three colorings of the same graph with different choices of CR-graph, indicated by bold lines. R, G, B and Y are colors. The multiple colors on some vertices indicate that these are coloring constrained but not color fixed. E.g. the left choice gives a 3-cycle that is constrained to RGB, RB and BG, but being an odd cycle, it is also *fixed* to RGB as a whole. Assume that the central vertices are labled 1 to 6, starting from the top. All three choices of CR graph show that {2, 4, 6} is a color identical set. The first and the second choice also show that {1, 3} is a color identical pair, but the last choice does not, which means that this choice of CR graph does not satisfy (ii) for *A* = {1} and *B* = {3}.

That every color fixation pair in a 4-chromatic graph induces a color fixation embrace, means that

      - only an odd cycle can color fix a vertex directly,

      - only an edge can color fix an edge directly, and

      - only a vertex can color fix an odd cycle directly.

Given a color fixed odd cycle *C*, a color fixed edge *e* and a color fixed vertex *v*, clearly,

      - a vertex that is adjacent to every vertex in *C*, is color fixed by *C*,

      - an edge whose vertices are adjacent to both vertices in *e*, is color fixed by *e*, and

      - an odd cycle whose every vertex is adjacent to *v*, is color fixed by *v*.

What is less clear, is if a direct color fixator *F* of a color fixee *E* in a 4-chromatic graph can contain *independent vertices*. However, any independent vertex in *F* would require its own color fixator *F'*, which must be color fixed in turn, but in the end, the color fixation of *E*, *F* and *F'* must all be referable to a single coloring reference graph, which in a 4-chromatic graph must be a triangle.

A *vertex* in a color fixed triangle *T* can color fix an *odd cycle*, an *edge* in *T* can color fix another *edge*, and *T as a whole* can color fix a *vertex*. Such color fixations open for new color fixations, but they also require additional edges, and the cycles induced by these edges may separate some of the new color fixees from some of the ones already in place.

**Lemma 1**. *Given a color fixator F in a 4-chromatic graph G, if F is not color fixed as a whole, any two independently color fixed vertices in F must be adjacent.*

*Proof*: Let *E* be subgraph of *G* color fixed directly by *F*, and let the triangle *R* = *abc* be a coloring reference graph for *E* and *F*, fixed to *red*, *green* and *blue*.

In order for every vertex in *F* to be adjacent to every vertex in *E*, no two vertices in *F* can be separated by a cycle, unless *E* itself is a cycle and *F* is a pair. However, if *E* is an odd cycle and an independently color fixed vertex in *F* is adjacent to every vertex in *E*, then *E* is color fixed by that vertex alone, and *F* cannot have any other vertices.

Suppose *F* contains two non-adjacent vertices *u* and *v*, such that *u* is fixed to *yellow* directly by *R*, and *v* is somehow fixed to *red*. Since *u* and *v* are not adjacent, there must be another *yellow*-fixed vertex *w* that contributes to fixing *v* to *red*. If *w* is color fixed directly by *R*, *u* and *w* must be on opposite sides of *R*, and *v* must be on the same side as *w*. If, alternatively, a vertex *x* is fixed to *red* by *bcu*, *w* can be fixed to *yellow* by *bcx*. This would place *u*, *v* and *w* on the same side of *R*, but then *v* would be separated from *u* by both *bcw* and *bcx*.

If *u* and *v* are not adjacent and none of them were fixed directly by *R*, more color fixation mediating vertices would be required, and at least one intermediate color fixation would have to occur between *u* and *v*, inducing a cycle that separated *u* and *v*. □

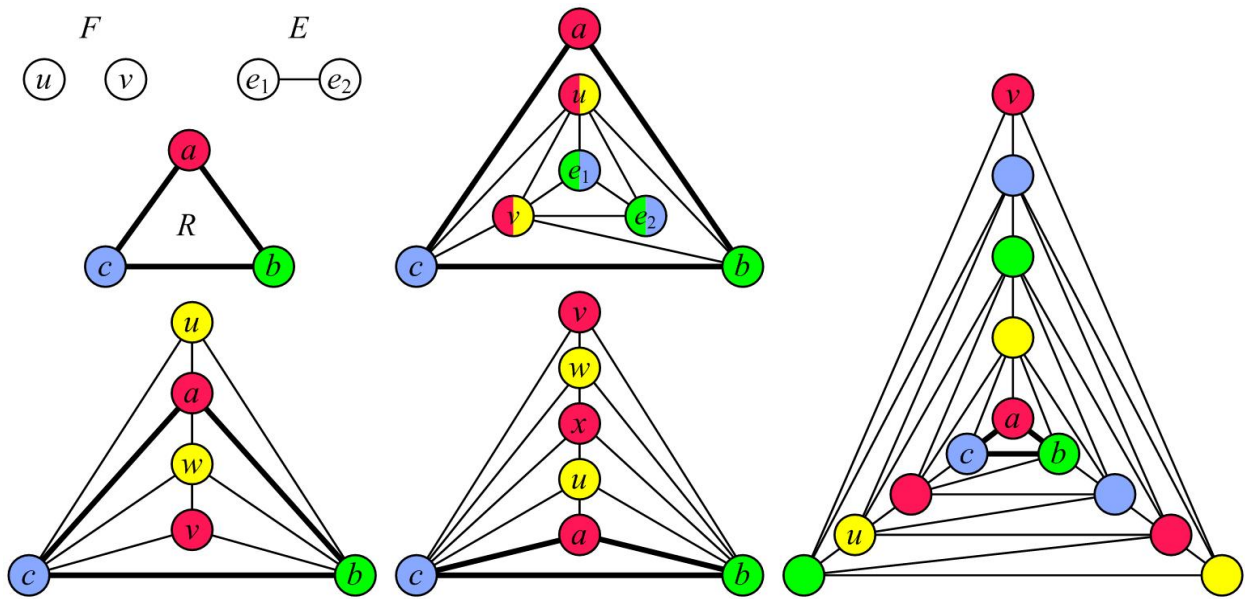

Figure 4. *R* = *abc* is the CR graph, *F* = {*u*, *v*} is the fixator, and *E* = $e_1 e_2$ is the fixee. If *u* and *v* are adjacent *F* can be fixed as a whole to {*red*, *yellow*} by the *green-blue* edge *bc* in *R*, and then *F* can color fix *E* to {*green*, *blue*}, but if not, additional vertices are required in order to fix *u* and *v* to different colors, and one way or the other *u* and *v* wind up being separated, making *F* unable to color fix *E*, no matter where *E* is placed.[4]

**Lemma 2**. *Only an odd cycle can color fix a vertex directly in a 4-chromatic graph.*

*Proof*: Let *v* be a vertex and let *F* be a direct color fixator of *v*.

If *F* is an odd cycle it can be color fixed as a whole and then color fix *v*.

If *V(F)* is a triple fixed independently to three different colors, then, by lemma 1, every two vertices in *F* must be adjacent, so we have an odd cycle.

If *F* contains more than three vertices and *F* is not an odd cycle, and thus not 3-critical, then, by lemma 1, any two non-adjacent vertices in *F* disqualifies *F* as a color fixator.  □

**Lemma 3**. *Only a vertex can color fix an odd cycle directly in a 4-chromatic graph.*

*Proof*: An odd cycle can be color fixed as a whole by a color fixed vertex.

Let $\mathcal{F} = \{A, B, C\}$ be a set of triangles color fixed by the triangle $R = abc$, and let $T = \{t, u, v\}$ be a set of vertices. The question is if *t*, *u* and *v* can be color fixed independently by $\mathcal{F}$, if *T* induces a triangle.

In order for *R* to color fix $\mathcal{F}$, *R* must be outside all of *A*, *B* and *C*, thus blocking any other vertices from being adjacent to every vertex in any of *A*, *B* or *C* from the outside. So, on the one hand, in order for *T* to be color fixed by $\mathcal{F}$, each of *t*, *u* and *v* must be inside one of *A*, *B* and *C* alone, and then *T* cannot induce a triangle. On the other hand, if *T* induces a triangle, *T* must also be outside all of *A*, *B* and *C*, but then none of *t*, *u* and *v* can be color fixed by either of *A*, *B* and *C*. The argument applies regardless of the number of triangles in $\mathcal{F}$, as long as there is at least one.

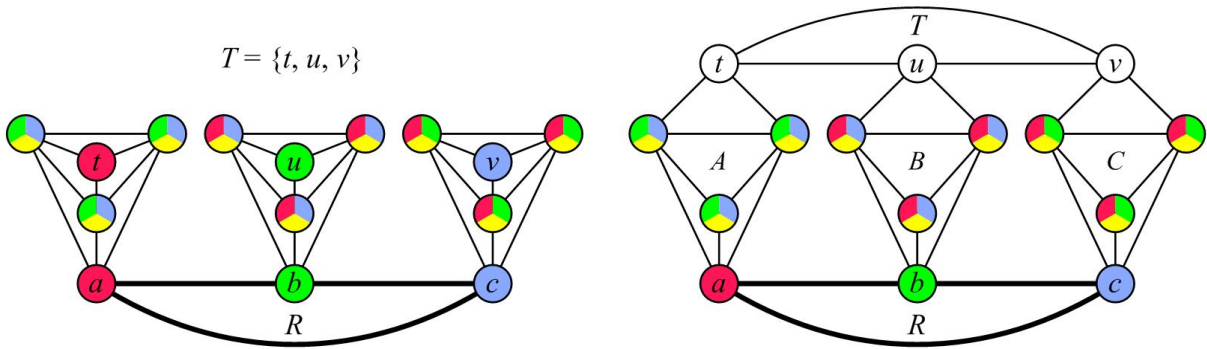

Figure 5.

Suppose instead that *T* is color fixed directly by *R*. If *a*, *b* and *c* are *red*, *green* and *blue*, and *t* is adjacent to all, *t* is fixed to *yellow*. That gives a *red*, *green* and *yellow* triangle *abt* that can fix *u* to *blue*, and that again gives a *green*, *blue* and *yellow* triangle *but* that can fix *v* to *red*—but instead of a direct color fixation of *T* by *R*, what we have, is a set of three direct color fixations of three different vertices in *T* by three different triangles in $R + T + \{at, au, bt, bu, bv, ct\}$.  □

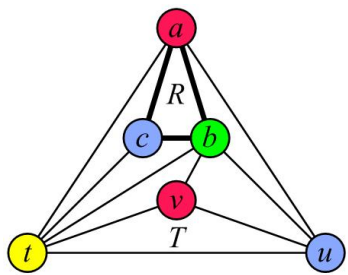

Figure 6.

**Definition 4**. 4-chromatic color fixation can propagate through a ***color fixation chain*** of alternating vertices and odd cycles that color fix each other—the ***vertex nodes*** and ***cycle nodes*** of the chain.[5] The minimal case comprises two vertices that are both adjacent to every vertex in the same triangle.

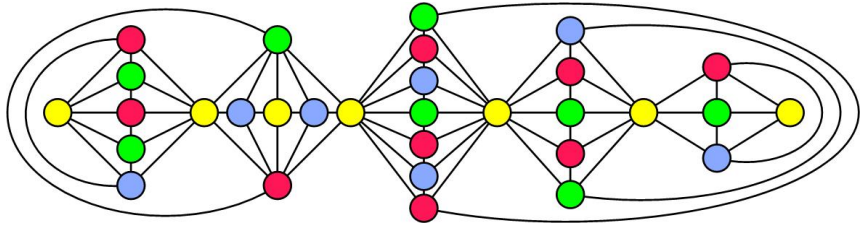

Figure 7. A 4-chromatic color fixation chain of vertices, here colored *yellow*, and odd cycles, here colored *red*, *green* and *blue*. Notice that the coloring of the cycle nodes is arbitrary, within the given constraints.

**Lemma 4**. *Two vertices in a 4-chromatic graph are color identical if and only if both are vertex nodes in the same color fixation chain.*

*Proof*: Since every pair of successive vertex nodes share an odd cycle as color fixator, and vice versa, the color fixation is repeated throughout the chain, and, since, by lemmas 2 and 3, only a vertex can color fix an odd cycle, and vice versa, color fixations cannot be transmitted between disjoint color fixation chains. □

**Definition 5**. Given two vertices $u$ and $v$ in a planar graph $G$, $u$ and $v$ are ***joinable*** if and only if they are not adjacent, and there is no cycle $X$ in $G$ such that $u$ and $v$ are on opposite sides of $X$. It follows that, if $uv$ is not an edge in $G$, then $G + uv$ is *planar* if and only if $u$ and $v$ are *joinable* in $G$.

**Lemma 5**. *There is no pair of joinable vertex nodes in a 4-chromatic color fixation chain.*

*Proof*: Given a color fixation chain $C$ of vertices and odd cycles, every cycle node $X$ partitions the vertex nodes in $C$ into two sets, the ones inside and the ones outside $X$, such that all the vertex nodes before $X$ in $C$ are on one side and all the vertex nodes after $X$ in $C$ are on the other side of $X$. □

**Theorem 1**: *There is no 4-chromatic planar graph with a joinable pair of color identical vertices.*

*Proof*: By lemma 4, color identical vertices must be in the same color fixation chain, and by lemma 5, there is no pair of joinable vertex nodes in a 4-chromatic color fixation chain. □

**Corollary 1:** *Every planar graph is 4-colorable.*

*Proof*: Suppose $G$ is a planar graph with chromatic number at least 5. We can assume that $G$ is 5-critical without loss of generality. Let $uv$ be an edge in $G$, and let $G' = G - uv$. Then $G'$ is 4-chromatic, and in every 4-coloring of $G'$, $u$ and $v$ must receive the same color, because otherwise we would also have a 4-coloring of $G$. Thus $u$ and $v$ are color identical and, by theorem 1, not joinable, contradicting the planarity of $G$.

In the other direction, suppose there is a 4-chromatic planar graph $G$ with a joinable $CI_4$ pair $\{u, v\}$. Then $G + uv$ is a 5-chromatic planar graph. □

## Notes

[1] My argument here is a variation of my argument in [1]. The parts of the argument that pertain to *separation* and *joinability* (called *adjaceability* in [1]) of vertices, have been simplified, the definitions of *color fixation* and of *coloring reference graph* have been sharpened, and in the proofs of the claims that *only a vertex can color fix an odd cycle*, and *vice versa*, the appeal to intuition has been dampened.

Notice also that the reference in one of these proofs to Fowler's proof that every uniquely 4-colorable graph is an Apollonian network [3], has been left out. One thing is that the reference is not needed, even if it is not irrelevant. Another thing is that Fowler's proof is computer assisted and, as such, of no use in what is meant to be a simple proof of an equivalent to the Four color theorem.

[2] My argument does not require a proof that *only an edge can color fix an edge directly in a 4-chromatic graph*, but to argue that *every 4-chromatic color fixation pair induces a color fixation embrace*, such a proof is needed, in addition to the proofs of lemmas 2 and 3:

*Proof*: Let $e$ be an edge and let $F$ be a direct color fixator of $e$. If $F$ is an edge, it requires two colors, and thus color fixes $e$—and even if the vertices of $F$ are independently fixed, by lemma 1, they must still be adjacent, thus $F$ must be an edge anyway.   □

[3] That every planar 4-chromatic graph contains a triangle, follows from the fact that the smallest triangle-free 4-chromatic graph, the Grötsch graph, is not planar [2]. (And as an à propos: it would seem that triangle free 4-chromatic graphs cannot contain color fixations.)

[4] In each of the three bottom graphs in Figure 4 every vertex is color fixed by a triangle, thus these graphs are Apollonian networks, and as such, uniquely 4-colorable. This is not relevant to my argument per se (see note 1), but it is related to it on a more general level, insofar as Fowler [3] has proven that every uniquely 4-colorable graph is an Apollonian network (without using that name), and shown that the Four color theorem is a corollary of this.

[5] There are also chains of edges that color fix edges, but such chains do not produce color identical pairs of vertices.